\documentclass [a4,12pt]{article}
\usepackage{wrapfig}
\usepackage[all]{xy}
\usepackage{amsmath,amssymb}
\usepackage{amsthm}
\usepackage{color}
\usepackage{url}
\usepackage{pb-diagram,pb-xy}

\author{Kohei Tanaka}

\title {\textbf{Reconstruction of Manifolds From Their Morse Functions}}

\newtheorem{theorem}{Theorem}[section]
\theoremstyle{definition}
\newtheorem{example}[theorem]{Example}
\newtheorem{definition}[theorem]{Definition}
\newtheorem{notation}[theorem]{Notation}

\theoremstyle{result}
\newtheorem{proposition}[theorem]{Proposition}
\newtheorem{lemma}[theorem]{Lemma}
\newtheorem{corollary}[theorem]{Corollary}

\newtheorem{main}{Main theorem}

\begin{document}

\newcommand{\rarrow}{\longrightarrow}
\newcommand{\larrow}{\longleftarrow}
\newcommand{\Mor}{\mathrm{Mor}}

\renewcommand{\thesection}{\arabic{section}}

\maketitle

A{\footnotesize BSTRACT}. 
This paper describes how to recover the topology of a closed manifold $M$ from a good Morse function $f$ on $M$.
The essential method was suggested by Cohen, Jones and Segal. They constructed a topological category $C_{f}$ 
and claimed that the classifying space $BC_{f}$ is homeomorphic to $M$.
We prove it from a different viewpoint with them using a cell decomposition of $M$ associated to $f$.
The cell complex $M_{f}$ equipped with the decomposition induces a topological category $C(M_{f})$ 
whose classifying space $BC(M_{f})$ is homeomorphic to $M$.
We show that $C(M_{f})$ is isomorphic to $C_{f}$ as a topological category. 
\vspace{0.3cm}






\section{Introduction}

It is well-known that important topological invariants of a 
closed smooth manifold $M$ can be obtained from a Morse function $f$ on $M$.

The Morse inequalities relate the Betti numbers $b_{i}(M)$ of $M$ with the number of critical points of $f$. As a consequence, 
the Euler characteristic $\chi(M)$ of $M$ is obtained by the sum 
$$\chi(M)=\sum_{i=0}^{\dim(M)}(-1)^{i}c(i)$$
of the number of critical points $c(i)$ with index $i$. 
These invariants $b_{i}(M)$ and $\chi(M)$ are originally based on the singular homology group $H_{*}(M)$ of $M$, 
and we can also obtain it from the Morse chain complex of $M$ generated by critical points of $f$.

Milnor gave the homotopy type of $M$ from $f$ in \cite{Mil}.
He constructed a cell complex which has the same homotopy type of $M$ 
and the same number of cells as the critical points of $f$.
Moreover, 
Cohen, Jones and Segal tried to recover the topology of $M$ from $f$ in \cite{CJS}. They constructed 
a topological category $C_{f}$ associated to $f$. 
It consists of 
critical points of $f$ as objects and the compactified moduli spaces of gradient flow lines of $f$ (see  Definition \ref{compact}) 
as morphisms.
They claimed that the classifying space $BC_{f}$ is homeomorphic to $M$ if 
$f$ satisfies the Morse-Smale condition.
However, their paper still remains to be a preprint. 
This paper gives a proof of the theorem 
by a viewpoint of cell decomposition of $M$ when $f$ is a bowl function (see Notation \ref{Morse}).

\begin{main}[Theorem \ref{ma}]
\textit{There is a homeomorphism $BC_{f} \cong M$ for a bowl function $f$ on $M$, 
where $BC_{f}$ is the classifying space of $C_{f}$.}
\end{main}

Dai Tamaki shows that a cell complex $X$ with a good structure called cylindrically normal can be recovered as the classifying space $BC(X)$
of a topological category $C(X)$ associated to $X$ in \cite{Tam}.
In order to apply his method to $M$, we give a cell decomposition of $M$ with a cylindrically normal structure 
for a bowl function $f$.
We call it the Morse theoretic decomposition of $M$ associated to $f$
and denote the cell complex $M$ with the decomposition by $M_{f}$.

In general, the Morse theoretic decomposition and its cylindrical structure depend on a function on the set of critical points of $f$.
It gives a gluing parameter and a gluing map in the sense that Qin described in \cite{Qin'}.  
Moreover, it induces an isomorphism between $C(M_{f})$ and $C_{f}$ as topological categories.

This paper is organized as follows. 

We recall classical Morse theory in section 2 including 
the cell complex $M_{f}$ based on Kalmbach's papers \cite{Kal}, \cite{Kal'}.

Section 3 describes the notion of cylindrically normal cell complex \cite{Tam} 
and shows the cell complex $M_{f}$ to be cylindrically normal.
It gives a topological category $C(M_{f})$ whose classifying space $BC(M_{f})$ is homeomorphic to $M$.

Section 4 gives a proof of the main theorem by showing that $C(M_{f})$ and $C_{f}$ are isomorphic as topological categories.
Since both categories have the same set of objects consists of critical points of $f$, 
we construct a homeomorphism between the spaces of morphisms $C_{f}(p,q) \cong C(M_{f})(p,q)$ 
compatible with composition for each critical points $p,q$.

\section{Morse theoretic decomposition of $M$}

Let $f$ be a Morse function on a closed Riemannian manifold $M$.
We give a cell decomposition of $M$ associated to $f$ following Kalmbach's paper \cite{Kal}.
It is obtained by the following steps.
\begin{enumerate}
\item We choose tubular neighborhoods of unstable manifolds of $f$,
\item define vector fields on the tubular neighborhoods and
\item construct characteristic maps from the flow lines of the vector fields.
\end{enumerate}

\begin{notation}\label{Morse}

We prepare some notations.

\begin{itemize}
\item $f$ is a Morse function on a closed smooth Riemannian manifold $M$ with dimension $m$.
\item $\mathrm{Cri}(f)$ is the set of critical points of $f$. For $p \in \mathrm{Cri}(f)$, let $\lambda(p)$ be the index and 
let $c_{p}$ be the critical value at $p$.
\item $-\Delta f$ is the negative gradient vector field of $f$. 
\item For a vector field $X$ on $M$ and $x \in M$, let $\gamma_{x}^{X}$ be the flow line of $X$ satisfying 
$\gamma_{x}^{X}(0)=x$. In particular, $\gamma_{x}$ denotes the gradient flow line $\gamma_{x}^{-\Delta f}$.
\item For a curve $\gamma : \mathbb{R} \rarrow M$, denote $s(\gamma)=\displaystyle\lim_{t \to -\infty}\gamma(t)$
and $e(\gamma)=\displaystyle\lim_{t \to +\infty}\gamma(t)$. 
\item $W^{s}(p)$ is the stable manifold and 
$W^{u}(p)$ is the unstable manifold at a critical point $p$.
For $p,q \in \mathrm{Cri}(f)$, define $p \geqq q$ by $W^{u}(q) \subset \overline{W^{u}(p)}$.
It gives a partial order of $\mathrm{Cri}(f)$. In particular, $p>q$ denotes $p \geqq q$ and $p \not=q$. 
\end{itemize}

For each critical point $p$,
there exists a neighborhood $U_{p}$ of $p$ in $M$ and a homeomorphism 
$\alpha_{p} : \mathbb{R}^{m} \rarrow U_{p}$ such that 
$$f_{p}(x,y)= -|x|^{2}+|y|^{2}+f(p)$$
and $\alpha_{p}(0)=p$ for $(x,y) \in \mathbb{R}^{\lambda(p)} \times \mathbb{R}^{m-\lambda(p)}=\mathbb{R}^{m}$ 
where $f_{p}= f \circ \alpha_{p}$.
We have
$$-\Delta f_{p}(x,y)= (2x,-2y), \ \gamma_{(x,y)}^{-\Delta f_{p}}(t)= \left(e^{2t}x,e^{-2t}y \right).$$

A Morse function $f$ is called a \emph{bowl function} if it satisfies the following two conditions for any $p,q \in \mathrm{Cri}(f)$.
\begin{enumerate}
\item The Morse-Smale condition, that is, $W^{s}(q)$ and $W^{u}(p)$ intersect transversally.
\item It satisfies either
$W^{u}(q) \cap \overline{W^{u}(p)} = \phi$ or $W^{u}(q) \subset \overline{W^{u}(p)}$.
\end{enumerate}
From now on we assume that $f$ is a bowl function through this paper. 

For $a>0$, define 
\begin{itemize}
\item $D^{m}(a)=\{x \in \mathbb{R}^{m}\ |\ |x| \leqq a\}$ and $D^{m}_{p}(a) = \alpha_{p}(D^{m}(a))$, 
\item $D^{\lambda(p)}_{p}(a) = \alpha_{p} \left(D^{\lambda(p)}(a) \times \{0\}\right)$ and $D^{m-\lambda(p)}_{p}(a) = \alpha_{p} \left(\{0\} \times D^{m-\lambda(p)}(a) \right)$.
\item $V_{p}^{u}(a) = \{ \gamma_{x}(t)\ |\ x \in \mathrm{Int}(D^{m}_{p}(a)) , t \in [0,\infty)\}$ and $V_{p}^{s}(a) = \{ \gamma_{x}(t)\ |\ x \in \mathrm{Int}(D^{m}_{p}(a)) , t \in (-\infty,0]\}$.
\end{itemize}
Take $\varepsilon_{p}>0$ for each critical point $p$ satisfying the following two conditions:
\begin{enumerate}
\item $D_{p}^{m}(\varepsilon_{p}) \cap D_{q}^{m}(\varepsilon_{q}) = \phi$ for $p \not= q$.
\item $V_{p}^{u}(\varepsilon_{p}) \cap V_{q}^{u}(\varepsilon_{q}) = \phi$ 
and $V_{p}^{s}(\varepsilon_{p}) \cap V_{q}^{s}(\varepsilon_{q}) = \phi$ 
for $\lambda(p)=\lambda(q)$.
\end{enumerate}
We write $\varepsilon = \min_{p \in \mathrm{Cri}(f)} \{\varepsilon_{p}\}$, $V_{p}^{u}=V_{p}^{u}(\varepsilon)$ and $V_{p}^{s}=V_{p}^{s}(\varepsilon)$.

\end{notation}

\begin{definition}
For a critical point $p$, define a diffeomorphism
$$\beta_{p}^{u} : \mathbb{R}^{\lambda(p)} \rarrow W^{u}(p)$$ 
by the following commutative diagram.
\[
\xymatrix{
\mathbb{R}^{\lambda(p)}  \ar[r]^{\beta_{p}^{u}} & W^{u}(p)  \\
\mathrm{Int}\left(D^{\lambda(p)}\right) \ar[r]_{\alpha_{p}}^{\cong} \ar[u]^{u}_{\cong}  & \mathrm{Int}\left(D_{p}^{\lambda(p)}\right) 
\ar[u]_{\gamma^{u}}^{\cong}
}
\]
where $u(x)=x/(1-|x|)$, $\gamma^{u}(x)= \gamma_{x}(s_{x})$ and $e^{2s_{x}}= 1/(1-|\alpha_{p}^{-1}(x)|$ respectively.
Similarly, we obtain a diffeomorphism $\beta_{p}^{s} : \mathbb{R}^{m-\lambda(p)} \rarrow W^{s}(p)$ by the following commutative diagram
\[
\xymatrix{
\mathbb{R}^{m-\lambda(p)}  \ar[r]^{\beta_{p}^{u}} & W^{s}(p)  \\
\mathrm{Int}\left(D^{m-\lambda(p)}\right) \ar[r]_{\alpha_{p}}^{\cong} \ar[u]^{u}_{\cong}  & \mathrm{Int}\left(D_{p}^{m-\lambda(p)}\right) 
\ar[u]_{\gamma^{s}}^{\cong}
}
\]
where $\gamma^{s}(x)=\gamma_{x}(s'_{x})$ and $e^{-2s'_{x}}=1/(1-|\alpha_{p}^{-1}(x)|)$.
\end{definition}

Given a Riemannian manifold $M$ and a submanifold $N$, 
it is well-known that $N$ has a tubular neighborhood in $M$.
However, we need to make a careful choice of particular tubular neighborhoods of stable and unstable manifolds of $f$.

\begin{definition}[Tubular neighborhood]\label{neigh}
Let $E \rarrow X$ be a vector bundle over a space $X$ with a Riemannian metric on $E$ 
and 
$$r : X \rarrow \mathbb{R}_{+}=\{x \in \mathbb{R}\ |\ x \geqq 0\}$$
be a continuous function. 
Define an open neighborhood 
$$U(E,r)=\{ v \in E_{x}\ |\  ||v|| < r(x), x \in X\}$$
of zero section $X_{0}$.
Let $N$ be a submanifold of $M$ and $\pi_{N} :  \nu(N) \rarrow N$ be 
the normal bundle of $N$.
If there exists a 
smooth function $r : N \rarrow \mathbb{R}_{+}$ and an embedding
$h : U(\nu(N),r) \rarrow M$ such that $h(N_{0})=N$, 
we call the image of $h$ 
a \emph{tubular neighborhood} of $N$.
Since stable and unstable manifolds are contractible, 
their normal bundles are trivial.
\end{definition}

\begin{definition}\label{defi}
Define $\rho_{n} : \mathbb{R}^{n} \rarrow \mathbb{R}_{+}$ by 
$$\rho_{n}(x)= \begin{cases} \sqrt{\varepsilon^{2}-|x|^{2}},  & \mathrm{if}\ 0 \leqq |x| \leqq \varepsilon /\sqrt{2}\\
\displaystyle\frac{\varepsilon^{2}}{2|x|}, & \mathrm{if}\ \varepsilon /\sqrt{2} \leqq |x| < \infty.
\end{cases} 
$$
This is a smooth function since
$$\frac{\partial}{\partial x_{i}} \sqrt{\varepsilon^{2}-|x|^{2}} = \frac{\partial}{\partial x_{i}} \displaystyle\frac{\varepsilon^{2}}{2|x|}  = -x_{i} \left(\displaystyle\frac{\varepsilon}{\sqrt{2}} \right)^{-1}$$
for $|x|= \varepsilon /\sqrt{2}$ and any $1 \leqq i \leqq n$.
Let $r_{u}$ denote the composition 
$$\rho_{\lambda(p)} \circ \left(\beta_{p}^{u}\right)^{-1} : W^{u}(p) \rarrow \mathbb{R}_{+}.$$
The open set 
$$V^{u} = \{ \gamma_{x}^{-\Delta f_{p}}(t) \ |\ x \in \mathrm{Int}(D^{m}(\varepsilon)), t \in [0,\infty) \}$$
is a neighborhood of $\mathbb{R}^{\lambda(p)}$ in $\mathbb{R}^{m}$.
For $(y,z) \in \partial V^{u} \subset \mathbb{R}^{\lambda(p)} \times \mathbb{R}^{m-\lambda(p)}$, we have
$\rho_{\lambda(p)}(y) = |z|$. 
Define a diffeomorphism $\xi^{u} : V^{u} \rarrow V^{u}_{p}$ by 
$\gamma_{x}^{-\Delta f_{p}}(t) \mapsto \gamma_{\alpha_{p}(x)}^{}(t)$. 
We obtain a homeomorphism 
$$h_{p}^{u} : U \left(W^{u}(p) \times \mathbb{R}^{m-\lambda(p)},r_{u} \right) \rarrow V_{p}^{u}$$
by the following commutative diagram
\[
\xymatrix{
U \left(W^{u}(p) \times \mathbb{R}^{m-\lambda(p)}, r_{u} \right) \ar[r]^{\hspace{1.5cm} h_{p}^{u}}& V_{p}^{u}  \\
U(\mathbb{R}^{\lambda(p)} \times \mathbb{R}^{m-\lambda(p)},\rho_{\lambda(p)})  \ar[u]^{\beta_{p}^{u} \times 1}_{\cong} \ar@{=}[r] & V^{u}, \ar[u]^{\cong}_{\xi^{u}}
}
\]
therefore $V^{u}_{p}$ is tubular neighborhood of $W^{u}(p)$. Similarly, we obtain a function on $W^{u}(s)$
$$r_{s}= \rho_{m-\lambda(p)} \circ \left(\beta_{p}^{s}\right)^{-1} : W^{s}(p) \rarrow \mathbb{R}_{+},$$ 
a neighborhood $V^{s}$ of $\mathbb{R}^{m-\lambda(p)}$ in $\mathbb{R}^{m}$ as
$$V^{s}=\{\gamma_{x}^{-\Delta f_{p}}(t) \ |\ x \in \mathrm{Int}(D^{m}(\varepsilon)), t \in (-\infty,0] \},$$
a diffeomorphism 
$\xi^{s} : V^{s} \rarrow V^{s}_{p}$ given by $\gamma_{x}^{-\Delta f_{p}}(t) \mapsto \gamma_{\alpha_{p}(x)}^{}(t)$, and 
a homeomorphism 
$$h_{p}^{s} : U \left(W^{s}(p) \times \mathbb{R}^{m-\lambda(p)},r_{s} \right) \rarrow V_{p}^{s}$$
by the following commutative diagram
\[
\xymatrix{
U \left(\mathbb{R}^{\lambda(p)} \times W^{s}(p), r_{s} \right) \ar[r]^{\hspace{1cm} h_{p}^{s}}& V_{p}^{s}  \\
U( \mathbb{R}^{\lambda(p)} \times \mathbb{R}^{m-\lambda(p)},\rho_{m-\lambda(p)})  \ar[u]^{1 \times \beta_{p}^{s}}_{\cong} \ar@{=}[r] & V^{s}. \ar[u]_{\cong}^{\xi^{s}}
}
\]
Hence $V_{p}^{s}$ is a tubular neighborhood of $W^{s}(p)$.
\end{definition}

\begin{lemma}\label{q}
For critical points $p>q$ (see Notation \ref{Morse}) of $f$, the intersection
$V_{q}^{s} \cap W^{u}(p)$ is a tubular neighborhood of $W(p,q)=W^{s}(q) \cap W^{u}(p)$ in $W^{u}(p)$.
\begin{proof}
The normal bundle of $W(p,q)$ in $W^{u}(p)$ is 
the pullback of the normal bundle of $W^{s}(q)$ along the inclusion $W(p,q) \hookrightarrow W^{s}(q)$.
Hence the normal bundle $\nu(W(p,q))$ is trivial since
so is the normal bundle $\nu(W^{s}(q))$.  
The restriction of the map $h_{q}^{s}$ in Definition \ref{defi} 
implies that $V_{q}^{s} \cap W^{u}(p)$ is a tubular neighborhood of $W(p,q)$ in $W^{u}(p)$.
\end{proof}
\end{lemma}

We define a vector field on $V^{u}_{p}$ and $M$.

\begin{definition}\label{Xp}
Recall that the diffeomorphism $\xi^{u} : V^{u} \rarrow V_{p}^{u}$ given in Definition \ref{defi}.
For a critical point $p$, define a vector field $X_{p}$ on $V_{p}^{u}$ by
$$X_{p}(x) =  \xi^{u}_{*}(0, -2x_{2}) \in T_{x}V_{p}^{u}$$
where $(x_{1},x_{2})= \left(\xi^{u} \right)^{-1}(x)$ and 
$\xi^{u}_{*} : T_{(x_{1},x_{2})}V^{u} \cong T_{x}V_{p}^{u}$.
\end{definition}

\begin{lemma}\label{p}
The vector filed $X_{p}$ satisfies the following conditions.
\begin{enumerate}
\item $X_{p}(f(x)) \leqq 0$, in particular $X_{p}(f(x))=0$ if and only if $x \in W^{u}(p)$.
\item $X_{p}(x)=0$ if and only if $x \in W^{u}(p)$.
\end{enumerate}
\begin{proof}
For $x \in V_{p}^{u}$, 
$$X_{p}(f(x)) =  - 4|x_{2}|^{2} \leqq 0$$
where $(x_{1},x_{2}) = \left(\xi^{u}\right)^{-1}(x) \in \mathbb{R}^{\lambda(p)} \times \mathbb{R}^{m-\lambda(p)}$.
It follows that $X_{p}$ satisfies the condition $1$. On the other hand, the definition of $X_{p}$ induces the condition $2$.
\end{proof}
\end{lemma}

\begin{definition}
Let $A, B \subset M$ be closed subspaces in $M$ with $A \cap B = \phi$.
A smooth function $g$ on $M$ is called a \emph{separation function} of $A$ and $B$ in $M$ if
$0 \leqq g(x) \leqq 1$, $g^{-1}(0)=A$, $g^{-1}(1)=B$ and $\Delta g(x) \not= 0$ for any $x \in M$ with $0< g(x) <1$.
\end{definition}

\begin{definition}\label{dis}
Let $\gamma$ be a curve on $M$ such that $\gamma(t) \cap D^{m}_{p}(\varepsilon) \not= \phi$ for some $t \in \mathbb{R}$.
Denote $|\gamma|_{p}$ as the distance between $\gamma$ and $p$ in $D^{m}_{p}(\varepsilon)$, 
that is, $|\gamma|_{p}= \min_{\gamma(t) \in D^{m}_{p}(\varepsilon)} \{ |\alpha_{p}^{-1}(\gamma(t))|\}$.
\end{definition}

\begin{lemma}\label{sep}
For $0<a<b<\varepsilon$,
there exists a separation function $g_{a,b}$ of $\overline{V^{u}_{p}(a)}$ and 
$V^{u}_{p}-V^{u}_{p}(b)$ in $V_{p}^{u}$.
\begin{proof}
Let $\ell$ be the function on $V_{p}^{u}$ given by $\ell(x)= |\gamma_{x}|_{p}$.
We have $\ell^{-1}([0, a))=V^{u}(a)$ and $\ell^{-1}([b, \varepsilon)) = V_{p}^{u}-V_{p}^{u}(b)$.
Define a smooth function
$$k(x) =\begin{cases} \mathrm{exp} \displaystyle\left(\frac{- 1}{x-a}\right),  & \mathrm{if}\ x>a\\ 
0, & \mathrm{if}\ x \leqq a
\end{cases}$$
for $x \in \mathbb{R}$.
The function $G(x) = \displaystyle\frac{k(x)}{k(x)+k(b+a-x)}$ satisfies 
$G(a)=0$, $G(b)=1$ and $0 \leqq G(x) \leqq 1$.
The composed map $G \circ \ell$ is our desired function $g_{a,b}$.
\end{proof}
\end{lemma}

Take a function $d$ on $\mathrm{Cri}(f)$ satisfying $d_{p}=d(p) \in (0,\varepsilon)$ for any $p \in \mathrm{Cri}(f)$.

\begin{proposition}\label{+}
For $n \geqq 0$, there exists a vector field $X_{n}^{d}$ on $V^{(n)} = \bigcup_{\lambda(p) \leqq n} V_{p}^{u}$ satisfying the following properties:
\begin{enumerate}
\item For $x \in D_{p}^{m}(\varepsilon)- \bigcup_{p<q} V^{u}_{q}(d_{q})$, we have $X_{n}^{d}(x)=X_{p}(x)$.
\item $X_{n}^{d}(f(x)) \leqq 0$, in particular, $X_{n}^{d}(f(x))=0$ if and only if $x \in W^{(n)}= \bigcup_{\lambda(p) \leqq n} W^{u}(p)$. 
\item $X_{n}^{d}(x)=0$ if and only if $x \in W^{(n)}$.
\end{enumerate}
\begin{proof}
We prove by induction on $n$.
When $n=0$, define $X_{0}^{d}=-\Delta f$. 
Assume that there exists a vector field $X_{n-1}^{d}$ on $V^{(n-1)}$ satisfying the properties.
For a critical point $p$ with index $n$, 
take the separation function $g_{p}=g_{d_{p}/2, d_{p}}$ on $V_{p}^{u}$ by Lemma \ref{sep}, 
and define a vector field $X_{n}^{d}$ on $V^{(n)}$ by
$$X_{n}^{d}(x) = g_{p}(x)X_{p}(x) + (1-g_{p}(x))X_{n-1}^{d}(x)$$
for $x \in V^{(n-1)} \cup V_{p}^{u}$.
Since $0 \leqq g_{p}(x) \leqq 1$, 
$$X_{n}^{d}(f(x)) = g_{p}(x)X_{p}(f(x)) + (1-g_{p}(x))X_{n-1}^{d}(f(x)) \leqq 0.$$
Lemma \ref{p} implies that $X_{n}^{d}$ has the desired properties.
\end{proof}
\end{proposition}

\begin{lemma}\label{vector}
For $n \geqq 0$, there exists a vector field $Y_{n}^{d}$ on $M$ satisfying the following properties:
\begin{enumerate}
\item For $p>q$, $\lambda(p)=n$ and $x \in W(p,q) - \bigcup_{p>r>q} V_{r}^{u}(d_{r})$, 
then $Y_{n}^{d}(x)=-\Delta f(x)$.
\item For a critical point $p$ with $\lambda(p) \leqq n$ and $x \in D^{m}_{p}(\varepsilon) - \bigcup_{p<q} V^{u}_{q}(d_{q})$, 
we have $(\alpha_{p}^{-1})_{*}(Y_{n}^{d}(x))= (2ax_{1},-2x_{2})$ for some $0 \leqq a \leqq 1$, where 
$(x_{1},x_{2}) = \alpha_{p}^{-1}(x) \in \mathbb{R}^{\lambda(p)} \times \mathbb{R}^{m-\lambda(p)}$.
\item $Y_{n}^{d}(f(x)) \leqq 0$, in particular, $Y_{n}^{d}(f(x))=0$ if and only if $x \in W^{(n-1)} \cup \mathrm{Cri}(f)$.
\item $Y_{n}^{d}(x)=0$ if and only if $x \in W^{(n-1)} \cup \mathrm{Cri}(f)$.
\end{enumerate}

\begin{proof}
The function $g = \prod_{\lambda(q) \leqq n-1} g_{q} : V^{(n-1)} \rarrow \mathbb{R}$ is a separation function such that
$g^{-1}(0) \supset W^{(n-1)}$ and $\overline{g^{-1}(1)} \supset \partial V^{(n-1)}$.
The vector field $Y_{n}^{d}$ on $M$ is given by 
$$Y_{n}^{d}(x)= g(x)\left(-\Delta f(x) \right)+ (1-g(x))X_{n-1}^{d}(x).$$
It satisfies the condition 2, 3 and 4 by Proposition \ref{+} and the property of $- \Delta f$.
If $x \in W(p,q) - \bigcup_{p>r} V_{r}^{u}(d_{r})$, then $Y_{\lambda(p)}^{d}(x)=-\Delta f (x)$ by the definition.
In the case of 
$$x \in  \left( W(p,q) - \bigcup_{p>r>q} V_{r}^{u}(d_{r}) \right) \cap D^{m}_{q}(d_{q}),$$
we have
$$(\alpha_{q}^{-1})_{*}(X_{q}(x)) =(0,-2x_{2}) = (\alpha_{q}^{-1})_{*}(-\Delta f(x)) \in T_{(x_{1},x_{2})}\left(\mathbb{R}^{\lambda(q)} \times \mathbb{R}^{m-\lambda(q)} \right)$$
where $(x_{1},x_{2}) = \alpha_{q}^{-1}(x) \in \mathbb{R}^{\lambda(q)} \times \mathbb{R}^{m-\lambda(q)}$.
Hence 
\begin{equation}
\begin{split}
Y_{\lambda(p)}^{d}(x)&= g(x)\left(-\Delta f(x) \right)+ (1-g(x))X_{\lambda(p)-1}^{d}(x) \\ \notag
&= g(x)\left((-\Delta f(x) \right)+ (1-g(x))\left( -\Delta f(x)\right) \\
&=-\Delta f(x).
\end{split}
\end{equation}
It implies that $Y_{n}^{d}$ satisfies the condition 1.
\end{proof}
\end{lemma}

\begin{notation}
For $p \in \mathrm{Cri}(f)$, $x \in M$ and $d \in (0,\varepsilon)^{\mathrm{Cri}(f)}$,
\begin{itemize}
\item $\gamma_{x,d}^{p}$ is the flow line of $Y_{\lambda(p)}^{d}$ satisfying $\gamma_{x,d}^{p}(0)=x$,
\item $Y^{u}_{d}(p)$ is the set $\{x \in M\ | \ s(\gamma_{x,d}^{p}) = p\}$ and 
\item $Y_{d}(p,q)$ is the set $\{x \in Y^{u}_{d}(p)\ | \ e(\gamma_{x,d}^{p})=q\}$.
\end{itemize}
Lemma \ref{vector} implies that $Y_{d}(p,q)=W(p,q) -  \bigcup_{p>r>q} V_{r}^{u}(d_{r})$.
\end{notation}

\begin{lemma}
For any critical point $p$ and $d \in (0,\varepsilon)^{\mathrm{Cri}(f)}$, we have $Y^{u}_{d}(p)=W^{u}(p)$.
\begin{proof}
Lemma \ref{vector} follows that $Y^{u}_{d}(p)-V^{(\lambda(p)-1)} = W^{u}(p)-V^{(\lambda(p)-1)}$.
If $x \in V^{(\lambda(p)-1)}$, we can find a critical point $q$ such that  $\lambda(q) < \lambda(p)$ and 
$x \in V^{u}_{q}- \bigcup_{r<q} V_{r}^{u}$.
Take $y,z \in D_{q}^{m}(\varepsilon)-V^{u}_{q}(d_{q})$ satisfying $\gamma_{x}^{p}(t)=y$ and 
$\gamma_{x}(t)=z$ for some $t \in \mathbb{R}$.
Lemma \ref{vector} implies that $y_{2}=z_{2}$ where $(y_{1},y_{2})=\alpha_{q}(y)$ and 
$(z_{1},z_{2})=\alpha_{q}(z)$ in $\mathbb{R}^{\lambda(q)} \times \mathbb{R}^{m-\lambda(q)}$.
If $x \in Y^{u}_{d}(p)$, then $y \in W^{u}(p) \cap V^{s}_{q}$. 
By Lemma \ref{q}, $W^{u}(p) \cap V^{s}_{q}$ is a tubular neighborhood of $W(p,q)$ in $W^{u}(p)$. It implies that  $z \in W^{u}(p) \cap V^{s}_{q}$ and 
 $x \in W^{u}(p)$.
Conversely, if $x \in W^{u}(p)$, then $z \in V_{q}^{s} \cap W^{u}(p)$. Therefore we have $y \in V_{q}^{s} \cap W^{u}(p)$ and $x \in Y^{u}_{d}(p)$.
\end{proof}
\end{lemma}

We define a cell decomposition of $M$ from the vector field $Y^{n}_{d}$.

\begin{definition}\label{chara}
Take a sufficiently small $\varepsilon'>0$ and identify $D^{\lambda(p)}$ and $D^{\lambda(p)}_{p}(\varepsilon') \subset M$.
A characteristic map 
$$\varphi_{p}^{d} : D^{\lambda(p)} \rarrow M$$
is given by $\varphi_{p}^{d}(x)=\gamma_{x,d}^{p}(t_{x})$ where $t_{x}=\tan \left(|x|-1/2 \right) \pi \in \mathbb{R}$. It gives a cell decomposition
$$M= \bigcup_{p \in \mathrm{Cri}(f)} Y_{d}^{u}(p) = \bigcup_{p \in \mathrm{Cri}(f)} W^{u}(p)$$
We call it the \emph{Morse theoretic decomposition} of $M$ associated to $f$ and 
denote the cell complex $M$ with the above decomposition by $M_{f}^{d}$.
Define a function $\epsilon$ on $\mathrm{Cri}(f)$ by $\epsilon_{p}=\varepsilon/2$ for any $p \in \mathrm{Cri}(f)$ 
and denote $M_{f}^{\epsilon}$ by $M_{f}$ simply.
\end{definition}

\section{A cylindrical structure on $M_{f}$}

Let $X=\bigcup_{\lambda \in \Lambda}e_{\lambda}$ be a cell complex. The face poset $P(X)$ of $X$
is $\Lambda$ with the partial order $\lambda \leqq \mu$ given by $e_{\lambda} \subset \overline{e_{\mu}}.$
We call a cell complex $X$ regular if any characteristic map 
is an embedding.
It is well-known that for a regular cell complex $X$, the classifying space (order complex) $BP(X)$ of the face poset is homeomorphic to $X$.
Unfortunately, 
the cell complex $M_{f}$ we gave in Section 2 is not regular in general.
Hence we construct a topological category $C(M_{f})$ whose classifying space $BC(M_{f})$ is homeomorphic to $M$, instead of poset.
In order to construct $C(M_{f})$, 
we need a good structure on $M_{f}$ called cylindrically normal in \cite{Tam}.

\begin{definition}[Stratified space]\label{def1}
Let $X$ be a topological space and $\Lambda$ be a poset. A \emph{stratification} of $X$ indexed
by $\Lambda$ is a surjective map $\pi : X \rarrow \Lambda$ satisfying the following properties:

\begin{enumerate}
\item Each $\pi^{-1}(\lambda)$ is connected and locally closed, i.e. it is open in $\overline{\pi^{-1}(\lambda)}$,
\item $\pi^{-1}(\lambda) \subset \overline{\pi^{-1}(\mu)}$ if and only if $\lambda \leqq \mu$.
\end{enumerate}
For simplicity, we denote $e_{\lambda} = \pi^{-1}(\lambda)$ and call it a \emph{stratum} with index $\lambda$.
Given a surjective map $\pi : X \rarrow \Lambda$, we have a decomposition of X, i.e.
\begin{enumerate}
\item $X= \bigcup_{\lambda \in \Lambda} e_{\lambda}$.
\item For $\lambda, \mu \in \Lambda$, $e_{\lambda} \cap e_{\mu}=\phi$ if $\lambda \not= \mu$.
\end{enumerate}

The indexing poset $\Lambda$ is called the face poset of $X$
and is denoted by $P(X)$. 
A \emph{stratified space} $(X,\pi)$ is a pair of space $X$ and its stratification $\pi$.

We say a stratum $e_{\lambda}$ is \emph{normal} if $e_{\mu} \subset \overline{e_{\lambda}}$ whenever
$e_{\mu} \cap \overline{e_{\lambda}} \not= \phi$.
When all strata are normal, the stratification is said to be normal.

Let $(X, \pi_{X})$ and $(Y,\pi_{Y})$ be stratified spaces. 
A morphism of stratified spaces $(f, \overline{f})$ is a pair 
of a continuous map $f : X \longrightarrow Y$ and a
map of posets $\overline{f} : P(X) \longrightarrow P(Y)$ making the following diagram commutative
\[
\xymatrix{
X \ar[d]_{\pi_{X}} \ar[r]^{f} & Y \ar[d]^{\pi_{Y}} \\
P(X) \ar[r]_{\overline{f}} & P(Y). }
\]
\end{definition}

\begin{example} The followings are some basic examples of stratified spaces.
\begin{enumerate}
\item A cell complex is a stratified space whose any stratum is homeomorphic to an open disk.
\item For any space $X$, we regard $\pi_{0}(X)$ as the trivial poset.
The projection $X \rarrow \pi_{0}(X)$ is a stratification on $X$.
We call it the trivial stratification on $X$. 
\item For two stratified spaces $(X,\pi_{X})$ and $(Y,\pi_{Y})$, the product $(X \times Y, \pi_{X} \times \pi_{Y})$ is a stratified space.
\end{enumerate}
\end{example}

\begin{definition}[Cylindrical structure]\label{cri}
A \emph{cylindrical structure} on a normal cell complex $X$ consists of
\begin{itemize}
\item a normal stratification on $\partial D^{n_{\lambda}}$ for each characteristic map 
$\varphi_{\lambda} : D^{n_{\lambda}} \rarrow X$,
\item a stratified space $P_{\lambda,\mu}$ and a morphism of stratified spaces
$$b_{\lambda, \mu} : P_{\lambda,\mu} \times D^{n_{\lambda}} \rarrow \partial D^{n_{\mu}} \subset D^{n_{\mu}}$$
for each $\lambda<\mu$ in $P(X)$, where we regard $D^{n_{\lambda}}$ as a trivial stratified space,
\item a morphism of stratified spaces 
$$c_{\lambda_{1},\lambda_{2},\lambda_{3}} : P_{\lambda_{2},\lambda_{3}} \times P_{\lambda_{1},\lambda_{2}} \rarrow P_{\lambda_{1},\lambda_{3}}$$
for a sequence $\lambda_{1}< \lambda_{2} < \lambda_{3}$ in $P(X)$
\end{itemize}
satisfying the following conditions:
\begin{enumerate}
\item The restriction of $b_{\lambda,\mu}$ to $P_{\lambda,\mu} \times \mathrm{Int}(D^{n_{\lambda}})$ is an embedding.
\item The following diagram is commutative
\[
\xymatrix{
P_{\lambda,\mu} \times D^{n_{\lambda}} \ar[d]_{b_{\lambda,\mu}} \ar[r]^{\ \ \ \mathrm{pr}_{2}} & D^{n_{\lambda}} \ar[d]^{\varphi_{\lambda}} \\
D^{n_{\mu}} \ar[r]_{\varphi_{\mu}} & X. }
\]
\item The following diagram is commutative
\[
\xymatrix{
P_{\lambda_{2},\lambda_{3}} \times P_{\lambda_{1},\lambda_{2}} \times D^{n_{\lambda_{1}}} 
\ar[r]^{\ \ \ 1 \times b_{\lambda_{1},\lambda_{2}}} \ar[d]_{c_{\lambda_{1},\lambda_{2},\lambda_{3}} \times 1 } 
& P_{\lambda_{2},\lambda_{3}} \times D^{n_{\lambda_{2}}} \ar[d]^{b_{\lambda_{2},\lambda_{3}}} \\
P_{\lambda_{1},\lambda_{3}} \times D^{n_{\lambda_{1}}} \ar[r]_{\ \ b_{\lambda_{1},\lambda_{3}}} & D^{n_{\lambda_{3}}}. }
\]
\item The map $c$ satisfies the associativity condition, i.e.
\[
\xymatrix{
P_{\lambda_{3},\lambda_{4}} \times P_{\lambda_{2},\lambda_{3}} \times P_{\lambda_{1},\lambda_{2}} 
\ar[r]^{\ \ \ 1 \times c_{\lambda_{1},\lambda_{2},\lambda_{3}}} \ar[d]_{c_{\lambda_{2},\lambda_{3},\lambda_{4}} \times 1} 
& P_{\lambda_{3},\lambda_{4}} \times P_{\lambda_{1},\lambda_{3}} \ar[d]^{c_{\lambda_{1},\lambda_{3},\lambda_{4}}} \\
P_{\lambda_{2},\lambda_{4}} \times P_{\lambda_{1},\lambda_{2}} \ar[r]_{\ \ c_{\lambda_{1},\lambda_{2},\lambda_{4}}} & P_{\lambda_{1},\lambda_{4}} 
}
\]
is a commutative diagram.
\item We have 
$$\partial D^{n_{\mu}} = \bigcup_{\lambda< \mu} b_{\lambda,\mu}(P_{\lambda,\mu} \times \mathrm{Int}(D^{n_{\lambda}}))$$
as a stratified space.
\end{enumerate}
A normal cell complex
equipped with a cylindrical structure is called a \emph{cylindrically normal} cell complex.
\end{definition}

Let's recall the definition of topological categories and functors between them.
In this paper, we treat only topological categories whose spaces of objects are discrete.

\begin{definition}
A topological category $C$ consists of the following data;
\begin{enumerate}
\item a set $C_{0}$ of objects and a space $C(x,y)$ of morphisms from $x$ to $y$ for any $x,y \in C_{0}$, 
\item a composition map $\circ : C(y,z) \times C(x,y) \rarrow C(x,z)$ for any $x,y,z \in C_{0}$,
\item an identity morphism $1_{x} \in C(x,x)$ for any $x \in C_{0}$,
\end{enumerate}
satisfying the identity condition and the associativity condition of composition.
A functor $F : C \rarrow D$
between two topological categories consists of a map
$F_{0} : C_{0} \rarrow D_{0}$ and continuous maps $C(x,y) \rarrow D(F_{0}x,F_{0}y)$ for each $x,y \in C_{0}$ 
preserving identity morphisms and composition.
The classifying space $BC$ of a topological category $C$ is a geometric realization of the nerve of $C$ (see \cite{Seg}).
\end{definition}

\begin{definition}[Cylindrical face category]
Let $X$ be a cylindrically normal cell complex. 
The \emph{cylindrical face category} $C(X)$ is defined by the following:
\begin{itemize}
\item The set of objects is $P(X)$.
\item The space of morphisms from $\lambda$ to $\mu$ is $P_{\lambda,\mu}$ and the composition is given by $c$.
\end{itemize}
\end{definition}

Tamaki showed the following theorem in \cite{Tam}.

\begin{theorem}[{[8, theorem 4.15]}]
\label{tam}
Let $X$ be a cylindrically normal cell complex, 
then the classifying space $BC(X)$ is homeomorphic to $X$.
\end{theorem}

The Morse theoretic decomposition is normal since $f$ is a bowl function.
In order to apply the above theorem to $M_{f}$, we show that $M_{f}$ has a cylindrical structure.

\begin{definition}
For $p,q \in \mathrm{Cri}(f)$, 
let $M(p,q)$ be the space of gradient flow lines of $f$ from $p$ to $q$ as a subspace of the mapping space 
$\mathrm{Map}(\mathbb{R},M)$ equipped with the compact open topology. It is called the moduli spaces of gradient flow lines from $p$ to $q$.
For any regular value $c_{q}<t<c_{p}$, there is a homeomorphism $W(p,q) \cap f^{-1}(t) \cong M(p,q)$ given by 
$x \mapsto \gamma_{x}$.
Let $B_{p}(a)$ be the space of gradient flow lines $\gamma$ satisfying $\mathrm{Im} \gamma \cap D_{p}^{m}(a) \not= \phi$.
For $d \in (0,\varepsilon)^{\mathrm{Cri}(p,q)}$, define
$$P_{d}(p,q)=M(p,q)- \bigcup_{p>r>q} \mathrm{Int}\left(B_{r}(d_{r}) \right).$$
In particular, we denote $P_{\epsilon}(p,q)$ by $P(p,q)$ simply, 
where $\epsilon$ is given in Definition \ref{chara}.
Lemma \ref{vector} induces $P_{d}(p,q) \cong \left(\varphi_{p}^{d}\right)^{-1}(q) \subset \partial D^{\lambda(p)}$, 
therefore $P_{d}(p,q)$ is a compact space.
\end{definition}

\begin{definition}\label{5}
For critical points $p>q$, define a vector field $Z_{d}$ on $\varphi_{p}^{-1}(W^{u}(q)) \subset \partial D^{\lambda(p)}$ by the following.
For $x \in \varphi_{p}^{-1}(W^{u}(q))$, take $t \in \mathbb{R}$ and $(y,z) \in W^{u}(q) \times \mathrm{Int} \left(D^{m-\lambda(q)}\right)$ 
satisfying 
$$h_{q}^{u}(y,z)=\gamma_{x,d}^{p}(t) \in \partial V^{u}_{q}(\varepsilon')$$
where 
$$h_{q}^{u} : W^{u}(q) \times \mathrm{Int}\left(D^{m-\lambda(q)} \right) \cong 
U(W^{u}(q) \times \mathbb{R}^{m-\lambda(q)}, r_{u}) \rarrow V^{u}_{q}$$
 is given in Definition \ref{defi}.
There exist open neighborhoods $U_{x}$ of $x$ in $\partial D_{p}$ and 
$V_{y} \times V_{z}$ of $(y,z)$ in $W^{u}(q) \times \partial D^{m-\lambda(q)}(|z|)$ 
and a diffeomorphism 
$\phi : U_{x} \stackrel{\cong}{\rarrow} V_{y} \times V_{z}$ such that $\phi(x)=(y,z)$.
Define $Z_{d}$ by $Z_{d}(x)=(\phi)^{-1}_{*}(Y_{\lambda(q)}^{d}(y),0)$ where 
$\phi_{*} : T_{x}U_{x} \rarrow T_{y}V_{y} \oplus T_{z}V_{z}$ is an isomorphism.
\end{definition}

\begin{lemma}\label{lift}
For $x \in W^{u}(q)$ and $\tilde{x} \in \partial D^{\lambda(p)}$ such that $\varphi_{p}^{d}(\tilde{x})=x$,
there exists a unique curve $\tilde{\gamma}_{x,d}^{p} : \mathbb{R} \rarrow \partial D^{\lambda(p)}$ such that 
$\varphi_{p}^{d} \circ \tilde{\gamma}_{x,d}^{p}=\gamma_{x,d}^{q}$ and 
$\tilde{\gamma}_{x,d}^{p}(0)=\tilde{x}$.
\begin{proof}
The desired curve $\tilde{\gamma}_{x,d}^{p}$ is obtained as the  
flow line $\gamma_{\tilde{x}}^{Z_{d}}: \mathbb{R} \rarrow \partial D^{\lambda(p)}$ of 
the vector field $Z_{d}$ in Definition \ref{5}.
\end{proof}
\end{lemma}

\begin{definition}\label{b}
Denote $\mathrm{Cri}^{p}=\{r \in \mathrm{Cri}(f)\ |\ p>r\}$, 
$\mathrm{Cri}_{q}=\{r \in \mathrm{Cri}(f)\ |\ r>q\}$ and $\mathrm{Cri}(p,q)=\mathrm{Cri}^{p} \cap \mathrm{Cri}_{q}$.
For $d \in (0,\varepsilon)^{\mathrm{Cri}(p,q)}$, $d' \in (0,\varepsilon)^{\mathrm{Cri}^{q}}$, 
and $t \in (0,\varepsilon)$, define $d*_{t}d' \in (0,\varepsilon)^{\mathrm{Cri}^{p}}$ by 
$$d*_{t}d'(i) = \begin{cases}
d(i), &\hspace{1cm} i \in \mathrm{Cri}(p,q)\\
t, &\hspace{1cm} i=q\\
d'(i), &\hspace{1cm} i \in \mathrm{Cri}^{q} \\
\varepsilon/2, &\hspace{1cm} \mathrm{otherwise}.
\end{cases}$$

Define a map
$$b_{t,d'} : P_{d}(p,q) \times D^{\lambda(q)} \rarrow  \partial D^{\lambda(p)}$$
for $t \in (0,\varepsilon)$ and $d' \in \mathrm{Cri}^{q}$ as follows.
For $(\delta,x) \in P_{d}(p,q) \times D^{\lambda(q)}$, 
take a sufficiently large $s \in \mathbb{R}$ and consider the point $(x,\delta(s)) \in D^{\lambda(q)} \times W(p,q)$.
Let $z$ denote $h_{p}^{s}(x,\delta(s)) \in V^{s}_{q} \cap W^{u}(p)=  V^{s}_{q} \cap Y^{u}_{d*d'}(p)$. 
Take the point $\tilde{x} \in \partial D^{\lambda(p)}$ satisfying $\tilde{x} = \gamma^{p}_{z,d*_{t}d'}(u)$ for some $u \in \mathbb{R}$. We have $\varphi_{p}^{d*_{t}d'}(\tilde{x})=x$.
Define $b_{t,d'}(\delta,x)= \tilde{\gamma}_{x,d*_{t}d'}^{p}(t_{x})$.
In particular, if $x \in \partial D^{\lambda(q)}$, then 
$b_{t,d'}(\delta,x) = e(\tilde{\gamma}_{x,d*_{t}d'}^{p}) \in \partial D^{\lambda(p)}$ and 
$\varphi_{p}^{d*_{t}d'}(b_{t,d'}(\delta,x)) = e(\gamma_{x,d*_{t}d'}^{q})$.
Furthermore, for $r<q<p$ and $t \in (0,\varepsilon)$ we define 
$$c_{t} : P_{d}(p,q) \times P_{d'}(q,r) \rarrow  P_{d*_{t}d'}(p,r)$$
by $c_{t}(\delta_{1},\delta_{2}) =  \gamma_{w,d*_{t}d'}^{p}$ where 
$w= b_{t,d'}(\delta_{1}, \delta_{2}(s))$ and 
$\delta_{2}(s) \in \partial D^{\lambda(q)}$.
We denote $c_{t}(\delta_{1},\delta_{2})$ by $\delta_{1} \circ_{t} \delta_{2}$.
\end{definition}

\begin{proposition}\label{compa}
For any $t, s \in (0,\varepsilon)$, $d_{1} \in (0,\varepsilon)^{\mathrm{Cri}(p,q)}$, 
$d_{2} \in (0,\varepsilon)^{\mathrm{Cri}(q,r)}$ and 
$d_{3} \in (0,\varepsilon)^{\mathrm{Cri}^{r}}$, the following diagram is commutative
\[
\xymatrix{
P_{d_{1}}(p,q) \times P_{d_{2}}(q,r) \times D^{\lambda(r)} 
\ar[r]^{\ \ \ 1 \times b_{s,d_{3}}} \ar[d]_{c_{t} \times 1 }
& P_{d_{1}}(p,q) \times D^{\lambda(q)} \ar[d]^{b_{t,d_{2}*_{s}d_{3}}} \\
P_{d_{1}*_{t}d_{2}}(p,r) \times D^{\lambda(r)} \ar[r]_{b_{s,d_{3}}} & D^{\lambda(p)}.}
\]
\begin{proof}
Let $d$ denote the function $d_{1}*_{t}d_{2}*_{s}d_{3}$.
For $(\delta_{1},\delta_{2},x) \in P_{d_{3}}(p,q) \times P_{d_{2}}(q,r) \times D^{\lambda(r)}$, 
$$b_{s,d_{3}} \circ (c_{t} \times 1)(\delta_{1},\delta_{2},x) = 
b_{s,d_{3}}(\delta_{1}\circ_{t} \delta_{2},x) = 
 \tilde{\gamma}_{x,d}^{p}(t_{x}) \in \partial D^{\lambda(p)}.$$
On the other hand, let $\eta$ be the curve on $\partial D^{\lambda(q)}$ given by 
$\eta(u) = \tilde{\gamma}_{x,d_{2}*_{s}d_{3}}^{q}(u)$
and $\tilde{\eta}$ be the curve on $\partial D^{\lambda(p)}$ given by
$\tilde{\eta}(v) = \lim_{u \to \infty} \tilde{\gamma}_{\eta(v),d}^{p}(u).$
It satisfies 
$$\varphi_{p}^{d} \circ \tilde{\eta}(v) =  \lim_{u \to \infty} \varphi_{p}^{d} \circ \tilde{\gamma}_{\eta(v),d}^{p}(u)
= \lim_{u \to \infty} \gamma_{\eta(v),d}^{q}(u)
= \varphi_{q}^{d}(\eta(v))
= \gamma_{x,d}^{r}(v).
$$
Both $\tilde{\eta}(0)$ and $\tilde{\gamma}_{x,d}^{p}(0)$ belong to 
$W^{u}(p) \cap V_{r}^{s}$ which is a tubular neighborhood of $W(p,r)$ in $W^{u}(p)$.
Let 
$\mathrm{pr}_{1} : W^{u}(p) \cap V_{r}^{s} \rarrow D^{\lambda(r)}$ and
$\mathrm{pr}_{2} : W^{u}(p) \cap V_{r}^{s} \rarrow W(p,r)$ 
be the projections, then
$$\mathrm{pr}_{1}(\tilde{\eta}(0)) =x= \mathrm{pr}_{1}( \tilde{\gamma}_{x}^{d}(0)),\  
\mathrm{pr}_{2}(\tilde{\eta}(0)) = \delta_{1}\circ_{t}\delta_{2}(u)= \mathrm{pr}_{2}( \tilde{\gamma}_{x}^{d}(0)) 
$$
where $\delta_{1}\circ_{t}\delta_{2}(u) \in \partial D_{p}$.
It follows that $\tilde{\eta}(0)= \tilde{\gamma}_{x,d}^{p}(0)$ and
$\tilde{\eta} = \tilde{\gamma}_{x,d}^{p}$.
Hence
$$b_{t,d_{1}*_{t}d_{2}} \circ (1 \times b_{s,d_{3}})(\delta_{1},\delta_{2},x)=\tilde{\eta}(t_{x}) =\tilde{\gamma}_{x,d}^{p}(t_{x})=b_{s,d_{3}} \circ (c_{t} \times 1)(\delta_{1},\delta_{2},x).$$
\end{proof}
\end{proposition}

The above proposition induces the following associativity condition of $c$.
Qin also considered the associativity of gluing of flow lines in \cite{Qin'}.

\begin{corollary}\label{as}
Take $t, s \in (0,\varepsilon)$ and $(\gamma_{1},\gamma_{2},\gamma_{3}) \in P_{d_{1}}(p_{1},p_{2}) \times P_{d_{2}}(p_{2},p_{3}) \times P_{d_{3}}(p_{3},p_{4})$. We have $(\gamma_{1}*_{t} \gamma_{2})*_{s} \gamma_{3} = \gamma_{1}*_{t} (\gamma_{2}*_{s} \gamma_{3})$, that is, the following diagram is commutative
\[
\xymatrix{
P_{d_{1}}(p_{1},p_{2}) \times P_{d_{2}}(p_{2},p_{3}) \times P_{d_{3}}(p_{3},p_{4})
\ar[r]^{\ \ \ 1 \times c_{s}} \ar[d]_{c_{t} \times 1 }
& P_{d_{1}}(p_{1},p_{2}) \times P_{d_{2}*_{s}d_{3}}(p_{2},p_{4}) \ar[d]^{c_{t}} \\
P_{d_{1}*_{t}d_{2}}(p_{1},p_{3}) \times P_{d_{3}}(p_{3},p_{4}) \ar[r]_{c_{s}} & P_{d_{1}*_{t}d_{2}*_{s}d_{3}}(p_{1},p_{4}).}
\]
\end{corollary}

\begin{theorem}\label{cylindrical}
The cell complex $M_{f}$ given in Definition \ref{chara} is cylindrically normal.
\begin{proof}
We regard $P(p,q)$ as a trivial stratified space and 
$$\partial D^{\lambda(p)}= \bigcup_{q<p} \left(\varphi_{p}\right)^{-1}(W^{u}(q))$$
is a normal stratification on $\partial D^{\lambda(p)}$ where $\varphi_{p}=\varphi_{p}^{\epsilon}$.
We verify that the maps $b_{q,p}=b_{\varepsilon/2,\epsilon}$ and $c_{r,q,p}=c_{\varepsilon/2}$ given in Definition \ref{b} satisfy the conditions of cylindrical structure in Definition \ref{cri} .
We have
$$\varphi_{p} \circ b_{q,p}(\delta,x) = \gamma^{q}_{x}(t_{y})=\varphi_{q}(x),$$
where $\gamma^{q}_{x}= \gamma^{q}_{x,\epsilon}$. Thus the following diagram is commutative
\[
\xymatrix{
P(p,q) \times D^{\lambda(q)} \ar[d]_{b_{q,p}} \ar[r]^{\ \ \ \ \ \mathrm{pr}_{2}} & D^{\lambda(q)} \ar[d]^{\varphi_{q}} \\
D^{\lambda(p)} \ar[r]_{\varphi_{p}} & M.}
\]
Proposition \ref{compa} gives the compatibility of $c$ and $b$, furthermore Corollary \ref{as} implies the associativity of $c$.
Take $\tilde{x} \in \varphi_{p}^{-1}(W^{u}(p))$ and let $x$ denote $\varphi_{p}(\tilde{x}) \in Y^{u}_{\epsilon}(q)$.
For $y=\left(\varphi_{q}\right)^{-1}(x) \in \mathrm{Int} \left(D^{\lambda(q)}\right)$
and $z =s \left(\tilde{\gamma}_{x}^{p} \right) \in \partial D^{\lambda(p)}$,
we can see $b_{q,p}(\gamma_{z}^{q},y)=\tilde{x}$. It follows that 
$$\partial D^{\lambda(p)}= \bigcup_{q<p} \left(\varphi_{p}\right)^{-1}(W^{u}(q))
=\bigcup_{q<p} b_{q,p} \left( P(p,q) \times \mathrm{Int}\left(D^{\lambda(p)}\right) \right).$$
\end{proof}
\end{theorem}

We consider the cylindrical face category $C(M_{f})$ of $M_{f}$. The set of objects is $\mathrm{Cri}(f)$
and the space of morphisms between two critical points $p,q$ is $P(p,q)$.
Theorem \ref{tam} induces the following corollary.

\begin{corollary}\label{a}
There is a homeomorphism $BC(M_{f}) \cong M$ for a bowl function $f$ on $M$.
\end{corollary}

\section{The comparison of $C(M_{f})$ and $C_{f}$}

We compare $C(M_{f})$ with the topological category $C_{f}$ 
given by Cohen, Jones and Segal in \cite{CJS}. 
The set of objects of $C_{f}$ consists of the critical points of $f$ as same as $C(M_{f})$.
On the other hand, the spaces of morphisms of $C_{f}$ is given by the compactified moduli spaces of gradient flow lines.

\begin{definition}[Compactified moduli space]\label{compact}
An ordered set 
$$I = \{p,p_{1},\cdots,p_{n-1},q\}$$
is a critical sequence if $p$, $q$ and $p_{i}$
$(i = 1, \cdots , n - 1)$ are critical points and $p > p_{1}> \cdots > p_{n-1}>q$. 
Let $M_{I}$ be the space of products of moduli spaces
$$M(p,p_{1}) \times \cdots \times M(p_{n-1},q)$$
and $P_{I}$ be the space of products $P(p,p_{1}) \times \cdots \times P(p_{n-1},q)$.
Denote the set of critical sequences from $p$ to $q$ by $S(p,q)$.
The compactified 
moduli space of gradient flow lines from $p$ to $q$ is given by the coproduct
$$\overline{M}(p,q) = \coprod_{I \in S(p,q)} M_{I}$$
as a set, and its topology is given by the following.
An element $\Gamma=(\gamma_{0},\cdots,\gamma_{n}) \in \overline{M}(p,q)$ gives the broken
curve $\overline{\Gamma} = \gamma_{0} * \cdots * \gamma_{n}$ on $M$ where $*$ is concatenation. 
Suppose the critical values of $f$ divide $[f(q), f(p)]$ 
into $\ell+1$ intervals $[c_{i},c_{i+1}]$ $(i=0,\cdots, \ell)$, 
where $c_{0}=c_{q}$ and $c_{\ell+1}=c_{p}$.
Choose a regular value $a_{i} \in (c_{i},c_{i+1})$.
The curve $\overline{\Gamma}$ intersects with $f^{-1}(a_{i})$ 
at exactly one point $x_{i}(\Gamma)$.
The evaluation map 
$$E : \overline{M}(p,q) \rarrow \prod_{0 \leqq i\leqq \ell} f^{-1}(a_{i})$$
is given by $E(\Gamma)= (x_{0}(\Gamma),\cdots, x_{\ell}(\Gamma))$.
Put the unique topology in $\overline{M}(p,q)$ such that the evaluation map 
$E$ is an embedding.
Note that the topology does not depend on the choice of $a_{i}$.
\end{definition}

\begin{theorem}[{[5, Theorem 7.4.]}]
The compactified moduli space of gradient flow lines $\overline{M}(p,q)$ satisfies the following properties.
\begin{enumerate}
\item $\overline{M}(p,q)$ is a compact manifold with boundaries whose interior is $M(p,q)$.
\item The inclusion $\overline{M}(p,r) \times \overline{M}(r,q) \hookrightarrow \overline{M}(p,q)$
is an embedding for any $q<r<p$.
\end{enumerate}
\end{theorem}

\begin{definition}\label{cf}
Define a topological category $C_{f}$ as follows.
\begin{enumerate}
\item The set of object is $\mathrm{Cri}(f)$.
\item The space of morphisms $C_{f}(p,q)$ is $\overline{M}(p,q)$ for $p,q \in \mathrm{Cri}(f)$.
\end{enumerate}
The composition is given by the inclusion $\overline{M}(p,r) \times \overline{M}(r,q) \hookrightarrow \overline{M}(p,q)$.
\end{definition}

We show that $C(M_{f})$ and $C_{f}$ are isomorphic as topological categories.
Since the sets of objects of both categories are $\mathrm{Cri}(f)$,
we construct a homeomorphism 
$$C_{f}(p,q) \rarrow C(M_{f})(p,q)$$
which is compatible with composition for $p,q \in \mathrm{Cri}(f)$.

\begin{lemma}\label{t}
Take critical points $p>q$ and $0<t<s <\varepsilon$. Let $\tau, \sigma$ be the functions on $\mathrm{Cri}(f)$ defined by 
$\tau(q)=t$, $\sigma(q)=s$ and $\tau(r)=\sigma(r)=\varepsilon/2$ for $r \not=q$.
Assume that $x \in W^{u}(q)$ and $\tilde{x}, \tilde{x}' \in \partial D^{\lambda(p)}$ 
satisfy $\varphi_{p}^{\tau}(\tilde{x})= \varphi_{p}^{\sigma}(\tilde{x}')=x$.
Then for any $u \in (-\infty,+\infty]$, 
there exists $v \in \mathbb{R}$ such that $\tilde{\gamma}_{x,\tau}^{p}(u)=\tilde{\gamma}_{x,\sigma}^{p}(v)$.

\begin{proof}
Let $y$ denote $\tilde{\gamma}_{x,\tau}^{p}(u) \in D^{\lambda(p)}$.
Since the curve $\gamma_{y,\tau}^{p}$ comes in contact with $D^{m}_{q}(t)$, 
there exists $z=(z_{1},z_{2}) \in \partial D^{m}_{q} \subset \mathbb{R}^{\lambda(q)} \times \mathbb{R}^{m-\lambda(q)}$ 
such that $\gamma_{y,\tau}^{p}(w)=z$ for some $w \in \mathbb{R}$.
On the other hand, $\gamma_{y,\sigma}^{p}(w)=(az_{1},z_{2})$ for some $0 \leqq a < 1$ by the property of $Y^{p}_{\sigma}$.
Therefore, the point $\varphi_{\sigma}^{p}(y)= e(\gamma_{y,\sigma}^{p})$ lies in the image of $\gamma_{x,\tau}^{q}=\gamma_{x,\sigma}^{q}$.
We have $\varphi_{p}^{\sigma}(y)=\gamma_{x,\sigma}^{q}(v)$ for some $v \in \mathbb{R}$ and 
$\tilde{\gamma}_{x}^{\tau}(u) = y = \tilde{\gamma}_{x}^{\sigma}(v)$.
\end{proof}
\end{lemma}

\begin{definition}\label{mu}
Take $\eta>0$ satisfying $\varepsilon/2<\eta<\varepsilon$. We define a map
$$\mu : P(p,r) \times P(r,q) \times (0,\eta] \rarrow M(p,q)$$
by $(\gamma,\delta,t) \mapsto \gamma \circ_{t}\delta.$
\end{definition}

\begin{lemma}
The map
$$\mu : P(p,r) \times P(r,q) \times (0,\eta] \rarrow M(p,q)$$
is injective.
\begin{proof}
Assume that $\gamma \circ_{t} \delta=\gamma' \circ_{s} \delta'$. 
Since $|\gamma \circ_{t} \delta|_{r}=t$ and $|\gamma' \circ_{s} \delta'|_{r}=s$ (see Definition \ref{dis}), 
we have $t=s$. Take $t<t'<\varepsilon$ and let $\tilde{x}= b_{t,\epsilon}(\gamma, \delta(u))$, 
$\tilde{x}'= b_{t',\epsilon}(\gamma, \delta(u))$ where $\delta(u) \in \partial D^{\lambda(q)}$.
By Lemma \ref{t}, 
$$\gamma \circ_{t} \delta=e \left(\tilde{\delta}^{\tau}_{\delta(u)} \right)= \tilde{\delta}^{\tau'}_{\delta(u)}(v)$$
for some $v \in \mathbb{R}$. Similarly, $\gamma' \circ_{t} \delta'= \tilde{\delta'}^{\tau'}_{\delta'(u')}(v')$.
Since $\tilde{\delta}^{\tau'}_{\delta(u)}$ and $\tilde{\delta'}^{\tau'}_{\delta'(u')}$ intersect each other, 
these are equal. Hence $\delta=\delta'$ and $\gamma = \gamma'$.
\end{proof}
\end{lemma}

\begin{lemma}\label{em}
The extended map 
$$\overline{\mu} : P(p,r) \times P(r,q) \times [0,\eta] \rarrow M(p,q) \cup \left(M(p,r) \times M(r,q)\right)$$
of $\mu$ in Definition \ref{mu} is defined by $\overline{\mu}(\gamma, \delta, 0)=\gamma \circ_{0}\delta= (\gamma, \delta)$.
This is an embedding when we regard $M(p,q) \cup \left(M(p,r) \times M(r,q)\right)$ as a subspace of $\overline{M}(p,q)$.

\begin{proof}
It is suffices to show that $\overline{\mu}$ is continuous at 
$$(\gamma,\delta,0) \in P(p,r) \times P(r,q) \times [0,\eta].$$
Take a sequence $\{(\gamma_{i},\delta_{i},t_{i})\}_{i \in \mathbb{N}}$ in $P(p,r) \times P(r,q) \times (0,\eta]$ 
converging to $(\gamma,\delta,0)$. 
We show that the sequence  
$\{(x_{i},y_{i}) = E(\gamma_{i} \circ_{ t_{i}} \delta_{i})\}$
 converge to $(x,y) = E(\gamma,\delta)$
in $f^{-1}(c_{r}-\varepsilon') \times f^{-1}(c_{r}+\varepsilon')$ for a sufficiently small $\varepsilon'>0$, where $E$ is the evaluation map in Definition \ref{compact}.
Since $\lim_{i \to \infty}(\gamma_{i},\delta_{i}) = (\gamma,\delta)$, we have $\lim_{i \to \infty}E(\gamma_{i},\delta_{i}) = (x,y)$.
The curve $\gamma_{i} \circ_{t_{i}} \delta_{i}$ passes through the end point of  the curve $\tilde{\delta}_{i}$ on 
$(\varphi_{p}^{\tau_{i}})^{-1}(W^{u}(p)) \cong M(p,q)-B_{r}(t_{i}/2)$. 
Since $\lim_{i \to \infty}t_{i} =0$,
$$\lim_{i \to \infty} e(\tilde{\delta}_{i}) = \lim_{i \to \infty} \tilde{x}_{i},$$
where $\tilde{x}_{i}=\tilde{\delta}_{i}(0)$. Therefore,
$$\lim_{i \to \infty}(x_{i},y_{i}) = \lim_{i \to \infty}(E(\gamma^{p}_{\tilde{x}_{i}})) = \lim_{i \to \infty}(E(\gamma_{i},\delta_{i}))=(x,y).$$
\end{proof}
\end{lemma}

\begin{definition}\label{e}
The associativity of the composition map $c$ in Corollary \ref{as} gives an embedding
$$e_{I} : P_{I} \times [0,\eta]^{|I|}  \rarrow \overline{M}(p,q)$$
by 
$$e_{I}(\gamma_{1},\cdots,\gamma_{n};t_{1},\cdots,t_{n-1})=\gamma_{1} \circ_{t_{1}}  \cdots \circ_{t_{n-1}} \gamma_{n}$$
A (broken) curve in the image of $e_{I}$ intersects $D_{r}^{m}(\eta)$ for any $r \in I$.
\end{definition}

\begin{theorem}\label{l}
For critical points $q<p$, there exists a homeomorphism
$$w_{q,p} : \overline{M}(p,q) \rarrow P(p,q)= M(p,q) - \left(\bigcup_{p>r>q}B_{r}(\varepsilon/2)\right)$$
making the following diagram commutative
\[
\xymatrix{
\overline{M}(p,q) \times \overline{M}(q,r)\ \  \ar@{_{(}->}[d] \ar[r]^{w_{q,p} \times w_{r,q}} &\ \ P(p,q) \times P(q,r) \ar[d]^{c_{\varepsilon/2}} \\
\overline{M}(p,r) \ar[r]_{w_{r,p}} & P(p,r). }
\]
\begin{proof}
For $\Gamma \in \overline{M}(p,q)$, 
let $I(\Gamma)$ be the critical sequence 
$$p=p_{0}>p_{1}> \cdots >p_{n-1}>p_{n}=q$$
satisfying  $\overline{\Gamma}(s) \cap  D_{p_{i}}^{m}(\eta) \not= \phi$ for some $s \in \mathbb{R}$, 
and $t_{i}$ denotes $|\overline{\Gamma}|_{p_{i}}$ for $0<i<n$. Note that 
$|\overline{\Gamma}|_{r}=0$ if $\Gamma=(\gamma_{1},\gamma_{2}) \in M(p,r) \times M(r,q)$.
There exists $\gamma_{i} \in P(p_{i},p_{i+1})$ such that
$$\Gamma = \gamma_{1} \circ_{t_{1}}  \cdots \circ_{t_{n-1}} \gamma_{n}$$ 
by Definition \ref{e}.
Define 
$$\iota : [0,\eta] \rarrow [\varepsilon/2,\eta]$$ 
by $\iota(t)= (1-\varepsilon/2\eta)t+\varepsilon/2$ and 
$$w_{q,p}(\Gamma) = \gamma_{1} \circ_{\iota(t_{1})} \cdots\circ_{\iota(t_{n-1})}  \gamma_{n}.$$
If $I(\Gamma)=p>q$, then $w_{q,p}(\Gamma)=\Gamma$.
This is a homeomorphism since $\iota$ is. 
Corollary \ref{as} induces the commutativity of the following diagram
\[
\xymatrix{
\overline{M}(p,q) \times \overline{M}(q,r)\ \  \ar@{_{(}->}[d] \ar[r]^{w_{q,p} \times w_{r,q}} &\ \ P(p,q) \times P(q,r) \ar[d]^{c_{\varepsilon/2}} \\
\overline{M}(p,r) \ar[r]_{w_{r,p}} & P(p,r). }
\]
\end{proof}
\end{theorem}

\begin{theorem}\label{ma}
There is a homeomorphism $BC_{f} \cong M$ for a bowl function $f$ on $M$, 
where $BC_{f}$ is the classifying space of $C_{f}$.
\begin{proof}
By theorem \ref{l}, $C_{f}$ and $C(M_{f})$ are isomorphic as topological categories.
Corollary \ref{a} implies that $BC_{f} \cong BC(M_{f}) \cong M$.
\end{proof}
\end{theorem}

\vspace{0.5cm}

\textit{Acknowledgements}: I would like to thank professor Dai Tamaki, Katsuhiko Kuribayashi and Keiichi Sakai for their useful suggestions on this work.
Their comments really helped me.

\vspace{0.5cm}

\begin{center}
Kohei Tanaka

Department of Mathematical Sciences

Shinshu University

Matsumoto

390-8621

Japan

E-mail: k-tanaka@math.shinshu-u.ac.jp
\end{center}

\end{document}